\catcode`\@=11
\ifx\macrosdouze@@loaded\UNDEFINED  \let\macrosdouze@@loaded=\null
\else \catcode`\@=12 \endinput\fi
\magnification=\magstep1
%\magnification=1000
%
%     Le petit carre de fin de demonstration
%

% definitions suivantes obsoletes : voir fin de fichier
\xdef \vrul@#1#2#3{%
      \vrule width#1\jot height#2\jot depth#3\jot \relax}
\def \carreblanc{\hbox{\jot=.2\p@
     \vrul@1{25}0\vrul@{24}10\kern-25\jot
     \vrul@{25}{25}{-24}\vrul@1{25}0}}
%%%%%%%%%%%%%%%%%
%%%% Ces macros activent les commandes
%%%% \ltimes et \rtimes pour les produits semidirects

\edef\arobasque{\the\catcode`\@}
\catcode`\@=11
\newfam\extrafam@
\font\extratext=msbm10
\font\extrascript=msbm10 at 7pt
\font\extrascriptscript=msbm10 at5pt

\textfont\extrafam@\extratext
\scriptfont\extrafam@\extrascript
\scriptscriptfont\extrafam@\extrascriptscript

\newfam\Bbbfam@
\font\Bbbtext=msam10
\font\Bbbscript=msam10 at 7pt
\font\Bbbscriptscript=msam10 at5pt

\textfont\Bbbfam@\Bbbtext
\scriptfont\Bbbfam@\Bbbscript
\scriptscriptfont\Bbbfam@\Bbbscriptscript

\def\newsymbol#1#2#3#4#5{\let\next@\relax
 \ifnum#2=\@ne\edef\next@{\the\Bbbfam@}\else
 \ifnum#2=\tw@\edef\next@{\the\extrafam@}\fi\fi
 \mathchardef#1="#3\next@#4#5}

\newsymbol\ltimes 226E
\newsymbol\rtimes 226F

%\catcode`\@=11
%    Quelques fontes
%
\font\sevenbf=cmbx7

%
%     Interligne simple ou double
%

%
\let\epsilon=\varepsilon

%
%     Gras de secours
%
\def\pmb#1{\setbox0=\hbox{#1}%
\hbox{\kern-.04em\copy0\kern-\wd0
\kern.08em\copy0\kern-\wd0
\kern-.02em\copy0\kern-\wd0
\kern-.02em\copy0\kern-\wd0
\kern-.02em\box0\kern-\wd0
\kern.02em}}
%
%     Sous-tilde
%
\def\undertilde#1{\setbox0=\hbox{$#1$}
\setbox1=\hbox to \wd0{$\hss\mathchar"0365\hss$}\ht1=0pt\dp1=0pt
\lower\dp0\vbox{\copy0\nointerlineskip\hbox{\lower8pt\copy1}}}
%
%quelques def
%texte dans les equations

%

%petites majuscules dans les biblios
\def\maj#1#2,{\rm #1\sevenrm #2\rm{}}
\def\Maj#1#2,{\bf #1\sevenbf #2\rm{}}
%lemmes avec petites majuscules
\outer\def\lemme#1#2 #3. #4\par{\medbreak
\noindent\maj{#1}{#2},\ #3.\enspace{\sl#4}\par
\ifdim\lastskip<\medskipamount\removelastskip\penalty55\medskip\fi}
%demonstrations

\def\Remark #1. {\noindent{\Maj REMARK,\ \bf #1. }}
%grandes fl\`eches avec un truc dessous

%lemmes avec petites majuscules grasses
\outer\def\Lemme#1#2 #3. #4\par{\medbreak
\noindent\Maj{#1}{#2},\ \bf #3.\rm\enspace{\sl#4}\par
\ifdim\lastskip<\medskipamount\removelastskip\penalty55\medskip\fi}
\outer\def\Theoreme #1. #2\par{\medbreak
\noindent\Maj{T}{H\'EOR\`EME},\ \bf #1.\rm\enspace{\sl#2}\par
\ifdim\lastskip<\medskipamount\removelastskip\penalty55\medskip\fi}
%
% d\'eriv\'ee normale:

\advance\vsize by 2 true cm
\voffset-1 true cm
\def\Eqalign#1{\null\,\vcenter{\openup\jot\m@th\ialign{
\strut\hfil$\displaystyle{##}$&$\displaystyle{{}##}$\hfil
&&\quad\strut\hfil$\displaystyle{##}$&$\displaystyle{{}##}$
\hfil\crcr#1\crcr}}\,}
%
%     Definitions polices msbm, msam ; symboles speciaux supplementaires
%
%\input Ensembles.tex
\catcode`\@=11
\magnification=\magstep1
%\magnification=1000
%
%     Le petit carr=E9 de fin de d=E9monstration
%

% d=E9finitions suivantes obsol=E8tes : voir fin de fichier
\xdef \vrul@#1#2#3{%
      \vrule width#1\jot height#2\jot depth#3\jot \relax}
\def \carreblanc{\hbox{\jot=.2\p@
     \vrul@1{25}0\vrul@{24}10\kern-25\jot
     \vrul@{25}{25}{-24}\vrul@1{25}0}}
%
%    Quelques fontes
%
\font\sevenbf=cmbx7

%
%     Interligne simple ou double
%
%\def \simpleinterligne{\normalbaselines
     \abovedisplayskip=12pt plus 3pt minus 9pt
     \abovedisplayshortskip=0pt plus 3pt
     \belowdisplayskip=12pt plus 3pt minus 9pt
     \belowdisplayshortskip=7pt plus 3pt minus 4pt
     \relax

\let\epsilon=\varepsilon

%
%     Gras de secours
%
\def\pmb#1{\setbox0=\hbox{#1}%
\hbox{\kern-.04em\copy0\kern-\wd0
\kern.08em\copy0\kern-\wd0
\kern-.02em\copy0\kern-\wd0
\kern-.02em\copy0\kern-\wd0
\kern-.02em\box0\kern-\wd0
\kern.02em}}
%
%     Sous-tilde
%
\def\undertilde#1{\setbox0=\hbox{$#1$}
\setbox1=\hbox to \wd0{$\hss\mathchar"0365\hss$}\ht1=0pt\dp1=0pt
\lower\dp0\vbox{\copy0\nointerlineskip\hbox{\lower8pt\copy1}}}
%
%quelques def
%texte dans les =E9quations

%

%petites majuscules dans les biblios
\def\maj#1#2,{\rm #1\sevenrm #2\rm{}}
\def\Maj#1#2,{\bf #1\sevenbf #2\rm{}}
%lemmes avec petites majuscules
\outer\def\lemme#1#2 #3. #4\par{\medbreak
\noindent\maj{#1}{#2},\ #3.\enspace{\sl#4}\par
\ifdim\lastskip<\medskipamount\removelastskip\penalty55\medskip\fi}
%d=E9monstrations

\def\Remark #1. {\noindent{\Maj REMARK,\ \bf #1. }}
%grandes fl\`eches avec un truc dessous

%lemmes avec petites majuscules grasses
\outer\def\Lemme#1#2 #3. #4\par{\medbreak
\noindent\Maj{#1}{#2},\ \bf #3.\rm\enspace{\sl#4}\par
\ifdim\lastskip<\medskipamount\removelastskip\penalty55\medskip\fi}
\outer\def\Theoreme #1. #2\par{\medbreak
\noindent\Maj{T}{H\'EOR\`EME},\ \bf #1.\rm\enspace{\sl#2}\par
\ifdim\lastskip<\medskipamount\removelastskip\penalty55\medskip\fi}
%
% d\'eriv\'ee normale:

\advance\vsize by 2 true cm
\voffset-1 true cm
\def\Eqalign#1{\null\,\vcenter{\openup\jot\m@th\ialign{
\strut\hfil$\displaystyle{##}$&$\displaystyle{{}##}$\hfil
&&\quad\strut\hfil$\displaystyle{##}$&$\displaystyle{{}##}$
\hfil\crcr#1\crcr}}\,}
%
%
%\input Ensembles.tex
%\mathchardef\square="1\msa@03
\def\carreblanc{\hbox{$\vrule height0.8ex width0.8ex$}}
\catcode`\@=12  % retablir

\def\liminf{\mathop{\underline{\rm lim}}}

\newfam\msbmfam
\font\tenmsbm=msbm10\textfont\msbmfam=\tenmsbm
\font\sevenmsbm=msbm7 \scriptfont\msbmfam=\sevenmsbm
\font\fivemsbm=msbm5 \scriptscriptfont\msbmfam=\fivemsbm
\mathchardef\R="7852
\mathchardef\N="784E
\mathchardef\L="784C
\mathchardef\M="784D
\mathchardef\C="7843
\mathchardef\Q="7851
\mathchardef\Z="785A
\mathchardef\P="7850
\mathchardef\T="7854
\mathchardef\I="7849
\mathchardef\H="7848

{\catcode`\@=11
%\global
%\mathchardef\semidirect="3\msb@6E
}

\def\semidirect{\mathrel{
\vrule width0.35pt height1.1ex
depth-0.05ex\kern-0.41ex\hbox{$\times$}}}

\vskip 1cm

\centerline{\bf On the width of lattice-free simplices}

\medskip

\centerline{\bf Jean-Michel KANTOR}

\vskip 1cm

\noindent {\bf   I Introduction}

Integral polytopes (see [Z ]for the basic definitions )
are of interest in combinatorics, linear programming,
algebraic geometry-toric varieties [D,O], number theory
[K-L.]. We study here lattice-free simplices, that is
simplices intersecting the lattice only at their
vertices.

A natural question is to measure the ``flatness" of
these polytopes, with respect to integral dual vectors.
This (arithmetical)  notion plays a crucial role :

  -in the classification (up to affine unimodular maps)
of lattice-free simplices in dimension $3$ (see
[O,MMM].

  -in the construction  of a polynomial-time algorithm
for integral linear programming ( flatness permits
 induction on the dimension, [K-L]).

 Unfortunately there were  no known examples (in any
dimension) of lattice-free polytopes with width bigger
than 2 .We prove here the following

{\bf Theorem:}

 Given any  positive number $\alpha$ strictly inferior
to  $1\over e$, for d large enough there exists a
lattice-free simplex of dimension d and width superior
to $\alpha d$.

\medskip

The proof is non-constructive and uses replacing the
search for lattice-free simplices in $\Z^d$ by the search
for "lattice-free lattices " containg $\Z^d$ ("turning
the problem inside out",see par.{ II}),specializing in
the next step to  lattices of a simple
kind,depending on a prime number p.The existence of
lattice-free simplices with big width is then deduced by
elementary computations,through a sufficient
inequality involving the dimension $d$,the width $k$
and the prime p  (see (14)).

The author thanks with pleasure H. Lenstra for crucial
suggestions,  V. Guillemin and I. Bernstein for
comments.

\bigskip

\noindent {\bf Notations}

\noindent ${\cal P}_d$: The set of integral polytopes in
$\R^d$; if $P$ is such a polytope,$P$ is a convex
compact set,the set Vert$(P)$ of vertices of $P$ is a
subset of
$\Z^d$.

\noindent ${\cal S}_d$: The set of integral simplices in
$\R^d$. In particular $\sigma_d$ will denote the
canonical simplex with vertices at the origin and
$$e_i = (0,\ldots,0,1,0,\ldots,0)  \; \hbox{--- $1$ at
the
$i$-th coordinate---}$$

\noindent $G_d$: the group of affine unimodular maps:
$$G_d = \Z^d \rtimes\ GL(d,\Z)$$
acts on $\R^d$ (preserving $\Z^d$), ${\cal P}_d$, and
${\cal S}_d$.

A $d$-lattice $M$ is a lattice with
$$\Z^d \subset M \subset {1 \over m} \Z^d \quad
\hbox{for some} \; m \in \N^\star.$$

\medskip

\noindent {\bf II. Turning simplices inside out.}

\smallskip

\noindent II.1. Let $\sigma$ be an integral simplex of
dimension $d$ in $\R^d$, and $L$ the sublattice of
$\Z^d$ it generates:
$$L = \left\{\sum^r_{i=1} m_i a_i, \quad a_i \in
Vert (\sigma)
 \quad m_i \in \Z\right\}.$$
We assume for simplicity one vertex
of $\sigma$ at the origin. The following is obvious:

\smallskip

\noindent {\bf Proposition 1.} {\sl i) There exists a
linear isomorphism
$$\varphi : \R^d \rightarrow \R^d \qquad \varphi(x) = y
= (y_j)_{j=1,\ldots,n}$$
such that
$$\varphi(\sigma) = \sigma_d.$$
It is unique up to permutation of the $y_i$'s, and
$$\varphi(L) = \Z^d, \quad \varphi(\Z^d) = M$$
where $M$ is a $d$-lattice.

ii) Conversely, given a $d$-lattice $M$, there
exists a linear isomorphism
$$\psi : \R^d \rightarrow \R^d$$
such that
$$\psi(M) = \Z^d,$$
and the image of $\sigma_d$ by $\psi$ is an integral
simplex $\sigma$ of dimension $d$ corresponding to $M$
as in i).}

All d-lattices generate $\R^d$ as vector space
over $\R$, and the Proposition is an easy consequence of
this.

\smallskip

\noindent {\sl Remarks.}

  -Because $\varphi$ is an
isomorphism, the following indices are equal
$$[\Z^d : L] = [M : \Z^d].$$

The  determinant of the lattice $L$ is classically
the volume of the parallelotope built on a basis of
$L$. If $\sigma$ generates $L$ as above,
$$\leqalignno{
&\hbox{vol} (\sigma) = {1 \over d!} \, \hbox{det} \,
L&(1)\cr
&\hbox{det} \, M = {1 \over d! \hbox{vol}
\,\sigma}.&(2)\cr }$$

\medskip
   -Proposition 1 has a straightforward extension   to
the case of two lattices  $  L $ and  $ M $ with
  $$L \subset M \subset {1\over p}L$$

\noindent II.2. Lattice-free simplices and
their width

Recall the following [K, K-L].

\smallskip

\noindent {\bf Definition 1.} An integral polytope $P$
in $\R^d$ is {\sl lattice-free} if
$$P \cap \Z^d = Vert(P) $$

\smallskip

\noindent {\bf Definition 2.} Given an integral non-zero
vector $u$ in $(\Z^d)^\star$, the $u$-width of the
polytope $P$ of ${\cal P}_d$ is defined by
$$w_u(P) = \max_{x,y \in P} \, < u, x - y >.$$

The {\sl width of} $P$ is
$$w(P) = \inf_{u \in (\Z^d)^\star \atop u \not =
0} w_u(P).\leqno(4)$$

Remark :The width is the minimal length of all
integral projections  ${u(P)}$ for non-zero $u$.

 II.3    Known results on the width of lattice-free
polytopes  in dimension $d$:

\smallskip

\noindent $d = 2$  :

\noindent Lattice-free simplices are all integral
triangles of area 1/2 ;they are equivalent to
$\sigma_2$.This is elementary.

\smallskip
\noindent  $ d = 3 $

\noindent Lattice-free polytopes have width one; in
the case of simplices, this result has various proofs
and applications (it is known sometimes as the
``terminal lemma", see [F, M-S, O, Wh]).

\smallskip
\noindent  $ d =4 $

 All lattice-free  simplices have at least one
basic facet (face with codimension one) [W] -this fact
is not true in higher dimensions.

Examples :
There exist some interesting examples:

-L.Schl\"afli 's polytopes,studied by Coxeter [C];

 -A recent example given by H.Scarf [private
communication ]:the simplex in dimension 5 with vertices
the origin,the first four vectors $e_i$ and for last
vertex ($23,39,31,43,57 $),has width 3.

  - We  have found  with the help of a
computer ,some examples of width 2,3 and 4 in dimension 4
and 5[F-K].

\smallskip

No other result seems to be known, apart from the
following asymptotic result:

\smallskip

\noindent {\bf Proposition  2.} {\sl There exists a
universal constant $C$ such that for any lattice-free
polytope of dimension $d$
$$w(P) \leq C d^2.\leqno(5)$$}

\medskip

\noindent {\sl Proof.} The ``Flatness Theorem" of [K-L]
asserts that there exists $C$ such that any  convex
compact set $K$ in $\R^d$ with
$$K \cap \Z^d=\phi$$
satisfies
$$w(K) \leq C\, d^2 \leqno(6)$$
where $w$ is defined as in I.2.

 If $P$ is any lattice -free polytope,take a
point
 $a$ in the relative  interior of $P$  and apply the
previous Flatness Theorem to the homothetic
$\widetilde{P}$ of
$P$ with respect to
$a$  and  fixed ratio $\alpha$ strictly less than one
 Then formula (4) shows that the width of P ,which is
proportional to the width of $\widetilde{P}$
,is also bounded by  a function of type (6).

Remark :Recent results  of [B] show that (5)
is true with  a right hand side proportional to
 d log d .

\noindent II.4. Turning the width inside out

Let us first define a new norm on $\R^d$ : If
$$\eqalign{
&x=(x_1\ldots,x_d) \in \R^d\cr
&\|x\|_\infty  = \sup_i |x_i|\cr
}$$

Define
$$\|x\| = \sup_i (0,x_i)-\inf_i
(0,x_i)\leqno(7)$$
It is the support function of  the following symmetric
convex  compact set in $(R^d)^\star$:
  $$ K = \sigma_d - \sigma_d$$
(see [O],p.182).

\noindent{\bf Lemma 1.}{\sl $\|\ \|$ is a norm, and
$$\|x\|_\infty \leq \|x\| \leq 2 \|x\|_\infty$$}

 From now on  $(\sigma,M)$ are as  in Proposition 1.
 We can identify the dual of the lattice $M$ with a
subgroup of $(\Z^d)^\star$:
$$\eqalign{
&\Z^d \subset M, \Longrightarrow M^\star \subset
(\Z^d)^\star\cr
&\xi \in M^\star : \xi =
(\xi_1,...\xi^d) \in
(\Z^d)^\star\cr
}$$

\medskip

\noindent {\bf Definition 3.} {\sl Let
$$w(M) = \inf_{\xi \in \M^\star\atop\xi \not =
0} \|\xi\|\leqno(8)$$}

Then we have
\smallskip

\noindent {\bf Proposition 3.}
$$w(\sigma) = w(M)$$

\medskip

\noindent {\sl Proof.} The isomorphism $\varphi$
changes $\sigma$ into  $\sigma_d$, linear forms $u$ on
$\Z^d$ into  linear forms on $M$, and
$$\eqalign{
w_\xi(\sigma_d) &= \sup_{x \in \sigma_d} (\xi,x) -
\inf_{y \in \sigma_d} < \xi, y >\cr
&= \sup_i(0,\xi_i) - \inf_i(0,\xi_i)\cr
&= \|\xi\|.\cr
}$$

\medskip

\noindent II.5. With notations as in II.1, we have
$$\sigma \cap \Z^d = Vert \, \sigma \Longleftrightarrow
M \cap \sigma_d = Vert(\sigma_d)$$
We can conclude this  part by asserting that the
existence of an integral  lattice-free simplex of
dimension
$d$, volume $v/d!$ and width at least $k$ is equivalent
with the existence of a  d-lattice $M$,containing
$\Z^d$, with
$$\left\{\eqalign{
&M \cap \sigma_d = Vert \, \sigma_d\cr
&w(M) \geq k\cr
&det(M) = {1\over v}.\cr
}\right.\leqno(9)$$

\vskip 1cm

\noindent III. {\bf In search of lattice-free simplices}
(asymptotically)

\smallskip

\noindent III.1.  We restrict our study to
$d$-lattices of type:
$$M(y) = \Z^d + \Z {1 \over p} y \qquad y \in
\Z^d \qquad M \not = \Z^d\leqno(10)$$
where $p$ is a prime number; clearly this lattice
depends only on the class of $y$ in $(\Z/p\Z)^d$.

\noindent {\bf Lemma 2
} The set of lattices $M$ (for a
fixed $p$) can be identified with the space of lines in
$(\Z/p\Z)^d$.

In particular the number of such lattices is

$$m(d,p) = {p^d - 1 \over p - 1}.\leqno(11)$$
Let $f(d,p)$ be the number of lattices $M$ such as
(10) satisfying
$$M \cap \check \sigma _d \not = \phi\leqno(12)$$
where
$$\check\sigma_d = \sigma_d \backslash
Vert(\sigma_d)\,.$$
(The lattice $M$ intersects $\sigma_d$ in other points
than the vertices).

\medskip

\noindent {\bf Lemma 3.}
$$f(d,p) \leq {(p + 1) ... (p+d) \over d!} - (d + 1).$$
\medskip

\noindent {\sl Proof.}
$$\left.\eqalign{
&x \in \M(y) \cap \check\sigma_d\cr
&x = z + {my \over p}\cr
}\right\} \Longrightarrow (m,p) = 1.$$
Writing $my\over p $as the sum of an integral vector
and a remainder we get
$$\eqalign{
&x = z + z' + {\tilde y \over p}, \quad 0 \leq
\widetilde{y_i} < p, \quad \widetilde{y_i} \in \N\cr
&x \in \check\sigma_d\cr
&\Longrightarrow z + z' = 0.\cr
  & x ={\tilde y\over p },\cr
  &\tilde y \in p\check\sigma_d
\cap \Z^d }$$
But the vectors $y,my,\tilde y$ define the same
line  in $(\Z/p\Z)^d $.
This shows that the number of lattices $M(y)$ satisfying
(12) is less than the number of points in $p
\check\sigma_d \cap \Z^d$, given by the right hand
side of the lemma [ E].

Let now $g(d,p,k)$ be the number of lattices $M(y)$ as
in (10) with
$$w(M(y)) \leq k.$$

\noindent {\bf Lemma 4.}

$$g(p,d,k) \leq 2 [(k+1)^{d+1} - k^{d+1}]
p^{d-2}.\leqno(13)$$
  Proof:

 The assumption on the lattice means the existence of a
vector
$\xi$ in
$\Z^d$
$$\eqalign{
&\xi \not = 0 \qquad y = (y_1,...y_d), \quad \xi =
(\xi_1,...\xi_d)\cr
&\sum \xi _i y_i \in p \Z\cr
&\|\xi\| \leq k \Longrightarrow \|\xi\|_\infty \leq
k.\cr
}$$

   The number of integral points  $\xi $  of  norm
less or equal to $ k $ is
$$    n  (k,d ) =(k+1 )^{d +1}  - k^{d +1}  $$
  [Proof :

 Let
$$ m = \inf_i(0,\xi_i)
        M =\sup_i(0,\xi_i)  $$

The possible values of $m$ are
$$ m = -k,\ldots\,-1,0$$
  a/For all values except $0$ one of the
$x_i$ has value $m$,and the others can take any value
between $m$ and $m+k$.
For each $m$ the
number of possibilities is equal to
$$S_1 =[{k+1}]^d -k^d $$

  b/When
$$ m =0$$
all $x_i's$ are non-negative,and the contribution is
  $$S_2 =[k+1]^d$$

  Adding up the contributions we get

$$ n (k,d)  =k [(k+1)^d -k^d] +(k+1)^d
             =(k+1)^{d +1}-k^{d+1}.]$$

Going back to the proof of Lemma 4,choose a vector $\xi$
with norm smaller than k ( strictly less than p):this
implies that the linear form defined by
$\xi$:
$$\hat\xi : (\Z/p\Z)^d \rightarrow \Z/p\Z$$
is surjective, and its kernel has $p^{d-1}$ elements;
the number of corresponding lattices is
$$r(p,d) = {p^{d-1} - 1 \over p - 1} \leq 2p^{d-2}.$$

We can choose at most $n(k,d)$ vectors $\xi$.

Hence ,by the mean value theorem ,
$$g(p,d,k) \leq 2 [(k+1)^{d+1} - k^{d+1}]
p^{d-2} \leq  2(d+1) (k+1)^d p^{d-2}.\leqno(13)$$
\medskip

\noindent III.3. From Lemmas 3 and 4
 we conclude that for large $d$ and $k$ the following
condition
 $$     2 (d+1){(k+1)}^d p^{d-2} + {(p+d)^d \over d!}
< p^{d-1}\leqno(14)$$
ensures the existence of a
lattice
$M(y)$ of width greater than $k$, dimension $d$, and
$$M(y) \subset {1 \over p} \Z^d.$$
The following is well-known :

\noindent {\bf Lemma 5.}

Given any sequence of numbers $( a_d)$going to
infinity,there exists an equivalent sequence
$(p_d)$ of prime numbers.

Proof :
Given  $\epsilon$ strictly positive we know from the
prime number theorem that for d  large enough  there
exists  a prime number $p_d$  in the interval  $ [(1
-\epsilon)   a_d,{(1+\epsilon)a_d}]$.This implies

$$\mid{ p_d -a_d}\mid  < \epsilon  a_d $$
for d large enough .

Choose now $\alpha $ arbitrary -we will  soon fix it
- and  a sequence
$(p_d)$ of primes with
$$p_d \sim \alpha d!$$
and let us find $\alpha$ and a sequence $ (k_d)$
such that
$$2( d+1) {(k_d+1)}^d p_d^{d-2} < {1 \over 2}
p_d^{d-1}\leqno(15)$$
$${(p_d+d)^d \over d!} < {1 \over 2}
p_d^{d-1}\leqno(16)$$

  These two conditions imply (14)

The condition (16) is satisfied for  large enough $d$ if
$$\alpha < {1 \over 2}.$$

Indeed
$$p_d +d\sim \alpha d! ;$$
since
$$\alpha < {1 \over 2}$$
(16) follows if we can show that

 $$ (1+d/p)^{d-1} \rightarrow 1  $$
$$d\rightarrow \infty)$$
 where
$$p=p_d$$
But
$$log(1+d/p)^{d-1} \leq (d-1)  d/p \sim  d^2/\alpha d!
\rightarrow 0$$

 Then (15) becomes
$$ k_d +1   < [{1 \over 4 ( d+1)} p_d]^{1\over d}$$

This last expression is equivalent ,because of
Stirling's formula , to
$d\over e
$.Hence if we choose any sequence of integral numbers
$(k_d)$with
$$k_d <\beta d $$
with
$$0 < \beta < {1 \over
e}\leqno(17)$$ then  (15 ) and (16 ) are satisfied
for large d .

\medskip

\noindent {\bf Theorem.} {\sl For any $\beta$ strictly
less than $1/e$, there exists for sufficiently large
$d$ a sequence of lattice-free simplices of dimension
$d$ and width $w_d$,}
$$w_d > \beta d.$$

  Defining
$$   w(d) = sup_\sigma w (\sigma)$$

supremum taken over all lattice-free simplices of
dimension d ,then the previous Theorem amounts to :

$$ \liminf_{d\rightarrow \infty}{w(d)\over d}  \geq
{1\over e}$$

{\bf Final Remark.}
The study above raises the hope of improving the
bounds on the maximal width, by  introducing
more general lattices generated by a finite
number of   rational vectors, and replacing the prime p
by powers  in  (10) (Note the study of general lattices
of such type in [Sh]). Unfortunately -and rather
mysteriously- our computations in these new cases give
the {\bf same} bounds.

\vskip 1cm

\noindent {\bf References}

\medskip

\item{[B]} Banaszczyk - Litvak - Pajor
A. Szarek    The flatness theorem,  the Gaussian
projection and analogues of the reverse Santalo
inequality for non-symmetric convex bodies

\item {[C]} Coxeter H.S.M.  The polytope $2_{21}$ whose
$27$ vertices correspond to the lines on the generic
cubic surface,American J.M.62;457-486,1940.

\item{[D]}   Dais D.,  Enumerative combinatorics of
invariants of certain complex threefolds with trivial
canonical bundle,   Dissertation, Bonn 1994.
\item {[E]}  L.Euler
Introductio in analysis infinitorum,Vol.1,Ch.16
(1748),253-257

\item{[F]}     Frumkin  H.,    Description of elementary
three-dimensional polyhedra. Conference on statistical
and discrete analysis, Alma Ata, 1981 (Russian).

\item{[F.K.]} Fermigier S. - Kantor J.M., Exemples de
grande \'epaisseur, Paris 1997.

 \item{[K]} Kantor J.M. Triangulations of integral
polytopes and Ehrhart polynomials,
 Betr\"age f\"ur Algebra und
geometrie, 1997.

\item{[K.L.]} Kannan R.  - Lovasz L.,
Covering minima and lattice-point free convex bodies,
Ann.of Mathematics, 128, 577-602, 1988.

 \item{[MMM]} Mori K. - Morrison D. - Morrison
I., On four dimensional terminal quotient
singularities, Mathematics of computation
51, 1988, N$^\circ$ 184, 769-786.

  \item{[M S]}  Morrison D. - Stevens G., Terminal
quotient singularities in dimension three and
four, Proceedings of the Amer. Math. Soc. 90,
1984, N$^\circ$ 1, 15-20.

\item{[O.]}  Oda T.,   Convex bodies and algebraic
geometry, Ergeb. Springer-Verlag, 1988.

  \item{[S]} Scarf H.E. Integral polyhedra in three
space Math. Oper. res. 10-403-438, 1985.

  \item {[Sh ]} Shimura G.Introduction to the
arithmetic theory of automorphic functions,Princeton
Univ.Press 1971

\item{[W.]}  Wessels,   Die S\"atze von White
..Diplomarbe\"\i t, Bochum 1989.

  \item{[Wh]}  White G.K., Lattice
tetrahedra, Canad. Journal of math. 16 (1961), 389-396.
    \item {[Z ]} Ziegler G. Lectures on
polytopes,Springer-Verlag ,GTM 152,1996
\bigskip

\hskip 7cm Jean-Michel Kantor

\hskip 7cm Centre de Math\'ematiques de Jussieu

\hskip 7cm Universit\'e Paris 7

\hskip 7cm Tour 46 5e \'etage Boite 247

\hskip 7cm 4, place Jussieu

\hskip 7cm F-75252 PARIS CEDEX 05

\bye